\documentclass[fleqn]{mat01}
\usepackage{times,mathtimy,amssymb,psfig}
\usepackage[all]{xy}
\begin{document}

\setcounter{page}{299}
\firstpage{299}

\def\ZZ{{\mathbb Z}}
\def\QQ{{\mathbb Q}}
\def\GG{{\bf G}}
\def\G{{\underline{\bf G}}}
\def\ga{{\Gamma}}
\def\CC{{\mathbb C}}
\def\v{{\varphi}}
\def\gm{{\Gamma}}
\def\gg{{\underline{\bf g}}}
\def\k{{bf k}}
\def\p{{\bf p}}

\newcommand{\x}{{\underline{x}}}
\newcommand{\fg}{\mathfrak g}
\newcommand{\fk}{\mathfrak k}
\newcommand{\fp}{\mathfrak p}

\def\claim{\trivlist\item[\hskip\labelsep{{\it Claim.}}]}
\def\theo{\trivlist\item[\hskip\labelsep{{\bf Theorem.}}]}
\def\examp{\trivlist\item[\hskip\labelsep{{\it Examples.}}]}
\def\consequence{\trivlist\item[\hskip\labelsep{{\it Consequence.}}]}
\def\question{\trivlist\item[\hskip\labelsep{{\it Question.}}]}
\def\lem{\trivlist\item[\hskip\labelsep{{\it Lemma.}}]}

\newtheorem{theore}{Theorem}
\renewcommand\thetheore{\arabic{theore}}
\newtheorem{definit}[theore]{\rm DEFINITION}

\newtheorem{theor}[theore]{\bf Theorem}
\newtheorem{propo}[theore]{\rm PROPOSITION}
\newtheorem{lemm}[theore]{Lemma}
\newtheorem{coro}[theore]{\rm COROLLARY}
\newtheorem{rem}[theore]{Remark}
\newtheorem{exampl}[theore]{Example}

\title{The congruence subgroup problem\footnote{This is essentially a
transcript of the plenary talk given at the Joint India--AMS Mathematics
Meeting held in December 2003 in Bangalore, India.}}

\markboth{M~S~Raghunathan}{The congruence subgroup problem}

\author{M~S~RAGHUNATHAN}

\address{School of Mathematics, Tata Institute of Fundamental Research,
Homi Bhabha Road, Colaba, Mumbai 400 005, India\\
\noindent E-mail: msr@math.tifr.res.in}

\volume{114}

\mon{November}

\parts{4}

\Date{MS received 22 July 2004}

\begin{abstract}
This is a short survey of the progress on the congruence subgroup
problem since the sixties when the first major results on the integral
unimodular groups appeared. It is aimed at the non-specialists and
avoids technical details.
\end{abstract}


\maketitle

\vspace{2pc}

\noindent The group $SL(2, \ZZ)$ of $2 \times 2$ integral matrices of
determinant 1 is a group that crops up in different contexts in
mathematics. Its structure is understood. The group has a natural family
of normal subgroups (of finite index). If $I \subset \ZZ$ is a proper
non-zero ideal, the subgroup $\{g \in SL(2, \ZZ)|g \equiv 1 (\bmod\,
I)\}$ is a subgroup of finite index which we will denote as $SL(2, I)$.

It is in fact the kernel of the natural homomorphism of $SL(2, \ZZ)$
into the finite group $SL(2, Z/I)$. Towards the end of the 19th century,
the question was raised if there were other examples of normal subgroups
of finite index; and Fricke--Klein exhibited such subgroups. It turns
out that there is a {\it surjective} homomorphism $\varphi\hbox{:}\
SL(2, \ZZ) \to A_5$ (alternating group on 5 symbols. Let $\Gamma = {\rm
kernel} \ \varphi$.\vspace{.1pc}

\begin{claim}
The group $\Gamma$ cannot contain a subgroup of the form $SL(2,
I)$.\vspace{.7pc}
\end{claim}

In what follows $(k)$ will denote $k\ZZ$. If \hbox{$I = n \ZZ, n > 0$} and
\hbox{$n = \prod_p p^{\alpha_p}$} is the prime factorization of \hbox{$n, SL(2,
\ZZ/n\ZZ) = \prod_p SL(2, \ZZ/(p^{\alpha_p}))$.} The natural map $SL(2,
\ZZ) \to S(2, \ZZ/n\ZZ)$ is surjective. So any simple quotient of $SL(2,
\ZZ)$ is a quotient of $SL(2, \ZZ/(p^{\alpha_p}))$ for some prime $p$.
Now the kernel of the map $SL(2, \ZZ/(p^{\alpha_p})) \to SL(2, \ZZ/(p))$
is a $p$-group and $SL(2, \ZZ/(p))/(\pm Id)$ is simple and non-abelian if
$p \not=2$ or $3$. We conclude that any simple quotient of $SL(2,
\ZZ)/I$ is of the form $SL(2, \ZZ/(p))/(\pm Id)$ for some prime $p>3$.

On the other hand, $A_5$ is not isomorphic to $SL(2, \ZZ/(p))$ for any
prime $p$.

The order of $SL(2, \ZZ/(p)) = (p-1) \cdot p \cdot (p+1)$. So $SL(2,
\ZZ/(p))\not\simeq A_5$ unless $p = 5$. The two sylow subgroup of $SL(2,
\ZZ/(p))$ is cyclic while that of $A_5$ is $\ZZ/2 \times \ZZ/2$. So
kernel $\varphi$ does not contain any $SL(2, I)$ with $I$ a proper
non-zero ideal.

This phenomenon raises the following question.  First, a definition.

A subgroup $\Gamma \subset SL(n, \ZZ)$ $(n$ integer $\geq 2$) is a {\it
congruence subgroup} iff there is a proper non-zero ideal $I \subset
\ZZ$ such that $\Gamma \supset SL(n, I) = \{g \in SL(n, \ZZ)| g
\equiv 1 (\bmod\,I)\}$.

Are there subgroups $\Gamma$ of finite index (note that $SL(n, I)$ has
finite index in $SL(n, \ZZ)$) which are not congruence subgroups?

We saw above that the answer is `yes' when $n = 2$.

In 1962, Bass--Lazard--Serre and independently Mennicke discovered that
$SL(2, \ZZ)$ is exceptional. They proved the following theorem.

\begin{theo}{\it
If $n >2${\rm ,} every subgroup of finite index in $SL (n, \ZZ)$ is a
congruence subgroup.}\vspace{.7pc}
\end{theo}

The problem can be generalized. One can pose it for other groups of
integral matrices such as symplectic ones or orthogonal ones (for a
quadratic form over $\ZZ$). One may also replace $\ZZ$ by integers in a
number field or $S$-integers for a {\it set $S$ of primes} including all
the archimedean primes.

We will now give a very general formulation. Apart from the fact that
many naturally arising examples fall within the ambit of this
formulation, the formulation suggests techniques for the attack that the
special cases may not suggest that readily. For the general formulation
we introduce the following notations: $k$ will be a number field; $V$, a
complete set of mutually inequivalent valuations of $k$; $\infty$, the
set of archimedean valuations; $S$, a subset of $V$ containing $\infty$.

For $v \in V \backslash \infty, k_v$ is the completion of $k$ at $v$,
${\cal O}_v = $ integers in $k_v$, ${\cal O}_S = \{ x \in k | x \in
{\cal O}_v \ {\rm for} \ v \not\in S\}$, the ring of $S$-integers in
$k$, ${\cal O} = {\cal O}_{\infty} = $ integers in $k$ (when $k = Q,
{\cal O} = \ZZ)$.

Next we recall the definition of a linear algebraic group defined over
$k$. We regard $k$ as a subfield of $\CC$. A linear algebraic group $G$
(defined) over $k$ is a subgroup of $GL(n, \CC)$ which is also the set
of zeros of a (finite) collection of functions on $GL(n, \CC)$ of the
form $P(g_{ij}, \det g^{-1})$ where $g = (g_{ij})_{1 \leq i, j \leq n} \in
GL(n, \CC)$ and $P$ is a polynomial in $(n^2 +1)$ variables with {\it
coefficients in} $k$. We will call $G$ a $k$-algebraic group or simply a
$k$-group. We denote by $G(k)$ the group $G \cap GL(n, k)$.

\begin{examp}$\left.\right.$
\begin{enumerate}
\renewcommand\labelenumi{\arabic{enumi}.}
\item $GL(n, \CC)$ is evidently one (over any $k$).

\item $SL(n, \CC) = \{g \in GL(n, \CC) | \det\ g = 1 \}$.

\item $D(n) = \{g \in GL(n, \CC) | g_{ij} = 0 \quad {\rm for} \ i
\not= j\}$.

\item The group of upper triangular matrices in $GL(n)$.

\item The group of upper triangular matrices with all diagonal entries
equal to 1.

\item Let $F$ be a symmetric non-singular $n \times n$ matrix over $k$
and
\begin{equation*}
\hskip -.5cm O(F) = \{g \in GL(n, \CC) | \ ^tg F g = F\}.
\end{equation*}
The orthogonal group of $F$ is a $k$-group.

\item $SO(F) = \{g \in O(F) | \det \  g = 1\}$.

\item In Example 6,  if one takes $n=2$ and $F =
\begin{pmatrix}
0 &1\\
1 &0
\end{pmatrix}$, then  $a \mapsto
\begin{pmatrix}
a &0\\
0 &a^{-1}
\end{pmatrix}$ gives an isomorphism of $D(1)$ on $SO(F)$.

\item Note that the algebraic group
\begin{equation*}
\hskip -.5cm \{g \in GL(n, \CC) | g_{ij} - \delta_{ij} = 0 \quad {\rm
for} \ i \not= 1, g_{11}=1\}
\end{equation*}
is isomorphic to $\CC^{n-1}$ -- again a group over any $k$.

\item Let $D$ be a division algebra over $k$ and $e_1, e_2, \ldots,
e_{d^2}$ a basis of $D$ over $k$. Then $e_1, e_2, \ldots, e_{d^2}$ is a
basis over $\CC$ of $D \otimes_k \CC$. Let $R_i$ denote the
multiplication by $e_i$ on the right
\begin{equation*}
\hskip -.5cm R_i\hbox{:}\ D\otimes_k \CC \to D \otimes_k \CC;
\end{equation*}
they are elements of $GL(d^2, \CC)$. Define $G = \{g \in GL(d^2, \CC)
| g R_i = R_i g$ for all $i, 1 \leq i \leq d^2\}$. This is an algebraic
group over $k$.
\end{enumerate}
\end{examp}

We now make  the definition.

\begin{definit}$\left.\right.$\vspace{.5pc}

\noindent {\rm A subgroup $\Gamma$ of $G(k)$ is a $S$-{\it congruence}
subgroup if it contains a subgroup of the form $G \cap GL(n, I)$ with $I$
a proper non-zero ideal in ${\cal O}_S$, as a subgroup of finite index.
Note that ${\cal O}_S/I$ is finite so that $G\cap GL(n, I)$ has finite
index in $G\cap GL(n, {\cal O}_S)$; thus $GL(n, I)\cap G$ is an
$S$-arithmetic subgroup (see Definition below).}
\end{definit}

\begin{definit}$\left.\right.$\vspace{.5pc}

\noindent {\rm A subgroup $\Gamma$ of $G(k)$ is a $S$-arithmetic
subgroup if for some (hence any) $S$-congruence subgroup $\Gamma'$ of
$G(k), \Gamma \cap \Gamma'$ has finite index in both $\Gamma$ and
$\Gamma'$.

We say that subgroups $H_1, H_2$ of a group $H$ are {\it commensurable}
iff $H_1 \cap H_2$ has finite index in both $H_1$ and $H_2$.

A $k$-morphism of a $k$-group $G\subset GL(n, \CC)$ into a $k$-group
$G'$ in $GL(n', \CC)$ is a group morphism $f\hbox{:}\ G \to G'$ such
that for every $(i', j'), 1 \leq i', j'\leq n'$, the $(i', j')$th entry
$f_{i'j'}(g)$ of $f(g)$ and det $f(g)^{-1}$ are polynomials with
coefficients in $k$ in the entries $g_{ij}, 1 \leq i,j \leq n$ of $g$ and
det $g^{-1}\hbox{:}\ f_{i' j'}(g) = P_{i'j'}(g_{ij}, \det\, g^{-1})
\det\, f(g)^{-1} = D(g_{ij}, \det\, g^{-1})$. Note that these are polynomials in
$(g_{ij} - \delta_{ij})$ and $(\det\, g^{-1} -1)$ as well.}
\end{definit}

\begin{lem}{\it
The inverse image in $G(k)$ of an $S$-congruence {\rm (}resp.
$S$-arithmetic{\rm )} subgroup of $G'(k)$ under a $k$-morphism
$f{\rm :}\ G \to G'$ is a $S$-congruence {\rm (}resp.
$S$-arithmetic{\rm )} subgroup of $G(k)$.}\vspace{.7pc}
\end{lem}

\noindent Note that $f(G(k)) \subset G(k')$.

Let
\begin{equation*}
f_{i'j'}(g) = Q_{ij} (g_{ij} - \delta_{ij}, \det\  g^{-1} -1)
\end{equation*}
and
\begin{equation*}
\det\ f(g^{-1}) = D'(g_{ij} - \delta_{ij}, \det\  g^{-1} -1).
\end{equation*}
Let
\begin{equation*}
\{c_i = a_i/b_i | i \in E, a_i\in {\cal O}_S, b_i  \in {\cal O}_S,
b_i \not= 0\}
\end{equation*}
be the collection of all the (non-zero) coefficients of $Q_{ij}$ and
$D'$. Let $J'$ be a proper non-zero ideal in ${\cal O}_S$ and $J =
(\prod_{i\in I} c_i)\cdot J'$. Then one sees immediately that
\begin{equation*}
f (G \cap GL(n, J)) \subset G' \cap GL(n, J').
\end{equation*}
Hence the lemma.

\begin{consequence}
The notions of $S$-arithmetic and $S$-congruence subgroups of $G(k)$
depend only on the $k$-isomorphism class of $G$ (not on the realisation
of $G$ as a $k$-subgroup of $GL(n, \CC))$. Then the congruence subgroup
problem is:
\end{consequence}

\begin{question}
Is every $S$-arithmetic subgroup of $G$ a $S$-congruence subgroup?
\end{question}

The case $k = Q$ and $S = \infty$ is itself sufficiently challenging; so
if you are not comfortable with the more general situation, you can make
the assumption $k = Q, S = \infty$. In this case ${\cal O}_S = \ZZ$. In
the other extreme case when $S = V, {\cal O}_S = k (= Q$ if $k = Q)$. We
saw above that the answer is, `No', for $k = Q, S = \infty$ \ and \ $ G
= SL(2, \CC)$, `Yes', for $k = Q, S = \infty$ \ and \ $G = SL(n, \CC),
n>2$.

Very substantial progress has been made on this general question. To
describe the progress, we first describe a way of measuring the failure
of an affirmative answer to the above question formulated by Serre. We
make the group $G(k)$ into a topological group in two different ways.
Let ${\cal A}_S$ (resp. ${\cal C}_S$) be the collection of all
$S$-arithmetic (resp. $S$-congruence) subgroups in $G(k)$ and ${\cal
J}_{a,S}$ (resp. ${\cal J}_{c, s})$ the topology of the unique structure
of a topological group on $G(k)$ for which ${\cal A}_S$ (resp. ${\cal
C}_S$) is a fundamental system of neighbourhoods of the identity. Let
$\widehat{G}_{a,S}$ (resp. $\widehat{G}_{c,S})$ be the completion of
$G(k)$ with respect to the natural (left-invariant) uniform structure.
Since ${\cal J}_{c,S}$ is weaker than ${\cal J}_{a,S}$, the identity map
as a map of $G(k)$ is uniformly continuous from the topology ${\cal
J}_{a, S}$ to the topology ${\cal J}_{c, S}$. Consequently the identity
map extends to a continuous homomorphism of $\widehat{G}_{a,S}$ on
$\widehat{G}_{c,S}$. We have a commutative diagram
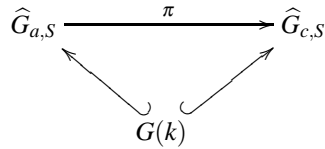
\begin{equation*}
\xymatrix{\widehat{G}_{a,S} \ar[rr]^-{\pi} & &\widehat{G}_{c, S}\\
&\ar@{_{(}->}[ul] G(k)\ \ar@{^{(}->}[ur] &\qquad . \\}
\end{equation*}
It turns out that $\pi$ is surjective and kernel $\pi := C(S, G)$ is
compact and totally disconnected. Evidently $C(S, G)$ provides a measure
of the failure of the family of $S$-arithmetic groups coinciding with
the family of $S$-congruence subgroups.

From the definitions it is not difficult to see that $C(S, G)$ is
contained in $\widehat{G}_a({\cal O}_S)$, the closure of $G({\cal O}_S)
= G \cap GL(n, {\cal O}_S)$ in $\widehat{G}_{a, S}(k)$ and is thus the
kernel of $\widehat{G}_a ({\cal O}_S) \to \widehat{G}_c ({\cal O}_S)$ $(= $
closure of $G({\cal O}_S)$ in $\widehat{G}_{c, S}(k))$. Because of this
one is able to conclude that $C(S, G)$ is totally disconnected and
compact: $\widehat{G}_a({\cal O}_S)$ is the profinite completion of
$G({\cal O}_S)$. We now pose the {\it congruence subgroup problem}:
\begin{equation*}
{\rm Determine} \ C(S, G) \tag{P}
\end{equation*}
for a given $G$ and $S$.

Observe that $C(S,G)$ is trivial iff every $S$-arith\-metic subgroup is
a $S$-congruence subgroup. Consider the extreme case $S = V$. Here
${\cal O}_S = k$; and $k$ has no proper non-zero ideals. So ${\cal O}_S
= \{G(k)\}$. Thus $C(S, G)$ is trivial if and only if $G(k)$ has no
proper (normal) subgroups of finite index. It is not difficult to see
that $G = GL(1)$ has lots of subgroups of finite index; so $C(V, G)$ in
general is non-trivial. However $C(V, G)$ has been conjectured to be
trivial under some natural restrictions on $G$ (Platonov--Margulis
conjecture).

The problem for general $G$ can be reduced to $G$ of a special kind
using the elaborate structure theory of linear algebraic $k$-groups.

Observe that $G \mapsto C(S, G)$ is a functor from the category of
$k$-groups into the category of compact totally disconnected ($\equiv$
profinite) groups.

\setcounter{theore}{0}
\begin{lemm}
If $G^o$ is the connected component of the identity in $G, G^o$ is a
$k$-group and the map $C(S, G^o) \to C(S, G)$ is an isomorphism.
\end{lemm}

\begin{lemm}
If a $k$-group $G$ is a semidirect product $(B\cdot N)$ with $B$ and $N
\ k$-subgroups and $N$ normal in $G$ and $C(S, N)$ is trivial{\rm ,}
then the map $C(S, B) \to C(S, G)$ is an isomorphism.
\end{lemm}

Lemma 2 combined with the structure theory of $k$-groups enables one to
reduce the problem to the case of {\it reductive} groups. A (connected)
$k$-group $G$ is {\it reductive} if it has no connected normal subgroups
consisting entirely of unipotent elements (unipotent $\equiv$ all
eigenvalues are $1$). It is a basic theorem that any $k$-group $G$ is the
semidirect product of a reductive $k$-group $B$ and the maximum normal
unipotent subgroup $R_uG$ (called the unipotent radical of $G$) which is
a $k$-group.

A unipotent $k$-group $U$ is a semidirect product $B\cdot U'$ where
$\dim U' = \dim U -1 $ and $B \simeq \hbox{Add} = \Bigg\{
\left. \begin{pmatrix}
1 &x\\
0 &1
\end{pmatrix} \right| x \in \CC \Bigg\} \subset GL(2, \CC)$.

Now $C(S, \text{Add})$ is trivial (an easy exercise).

An induction on dimension shows that $C(S,G)$ is trivial if $G$ is
unipotent (i.e. consists entirely of unipotent matrices). Lemma 2 thus
reduces our problem to the case of reductive $G$.

A $k$-group $T$ is a {\it torus} if it is connected and can be
conjugated into diagonal matrices in $GL(n, \CC)$.

It is again a basic result that if $G$ is a reductive $k$-group, $G$
contains a central $k$-torus $T$ such that $G/T$ has no connected
abelian normal subgroups. Information on $T$ and $G/T$ separately can be
pieced together to obtain results on $G$: this is somewhat delicate
though. We will now deal with $k$-tori. These are abelain but they need
much more subtle handling than unipotent groups. One has the following:

\begin{theo}\hskip -.5pc {\bf (Chevalley).}\ \ {\it $C(S, G) = \{1\}$
if $S$ is {\it finite} and $G$ is a torus.}\vspace{.7pc}
\end{theo}

This is false if $S$ is infinite. However one knows the structure of $T$
in sufficient detail to get considerable information on $C(S, G)$ in
this case too. Chevalley's theorem as also other information on $C(S,
G)$ for $S$ infinite needs some class field theory.

I will now say something on the most important case: $G$ semisimple,
i.e. $G$ has no non-trivial connected abelian normal subgroups. We will
make two more assumptions viz. that $G$ is simply connected and that $G$
is absolutely almost simple -- the latter means that $G$ has no proper
connected normal subgroups. This last assumption is not really
restrictive but simple connectivity is. However information in the
simply connected case can be effectively used to handle the general
case. To explain what is expected to be true, I need to introduce some
other concepts.

A $k$-split torus $T$ is $k$-torus $k$-isomorphic to the group $D(n)$
of all diagonal matrices in $GL(n)$ for some $n$.

It is a theorem of Borel--Tits that all maximal $k$-split tori in a
$k$-group $G$ are mutual conjugates under $G(k)$ and their common
dimension is called the $k$-rank of $G$. It is again a theorem of
Borel--Tits that $k$-rank $G \geq 1$ iff $G(k)$ has non-trivial unipotent\break
elements.

The $S$-rank of $G$ is the number $\sum_{v \in S} \ k_v$-rank $G$.

One expects the following: Assume $S$ is such that $k_v$-rank $G > 0$
for all $v \in S\backslash \infty$ and $S$-rank $G \geq 2$. Then $C(S,
G)$ is trivial or isomorphic to the group $\mu_k$ of roots of unity\break in
$k$.

Note that $k_v$-rank and $S$-rank $SL(n)\geq 2$ for any $v$ and any $S$
for $n \geq 3$. It is now known that the expectation is indeed true for
any $G$ with $k$-rank $G \geq 1$ (Bass, Lazard, Milnor, Serre, Mennicke,
Matsumoto, Deodhar, Vaserstein, Bak, Rehman, Prasad and Raghunathan).
The strategy in all this work has been to break the proof into two parts.

\begin{enumerate}
\item[(1)] Show that
\begin{equation*}
\hskip -.5cm 1 \to C(S, G) \to \widehat{G}_{a,S}(k) \to \widehat{G}_{c,
S}(k)\to 1
\end{equation*}
is a central extension.
\end{enumerate}

This has the consequence that $C(S, G)$ is the Pontrjagin dual of the
kernel of the map
\begin{equation*}
H^2(\widehat{G}_{c, S}(k), \QQ/\ZZ) \to H^2 (G(k), \QQ/\ZZ).
\end{equation*}
The first group is the cohomology group based on continuous co-chains,
the second is the usual group cohomology.

\begin{enumerate}
\item[(2)] Show that the above kernel is $\mu_k$ or trivial.
\end{enumerate}

The latter programme has in fact been carried out for all $G$ (Moore,
Matsumoto, Deodhar, Prasad, Raghunathan and Rapinchuk).

The centrality of the sequence in (1) above has also been settled in
many cases of $k$-rank 0 (groups of type $B_n, C_n, D_n$ and some
exceptional groups) (Kneser, Rapinchuk and Tomanov). The case $S = V$
for anisotropic $G$ is much more delicate than for isotropic $G$. The
first results here are due to Kneser. Anisotropic groups of type $A_1$
were dealt with by Platonov, Rapinchuk and Margulis. Groups of type
$B_n, C_n, D_n$ and some exceptional groups have been dealt with
(\^Cernusov, Rapinchuk, Sury and Tomanov). Groups of type $A_n$ pose the
greatest challenge. For {\it inner} forms of $A_n$ and $S = V$, $C(V,
G)$ has been determined. One can reformulate this as follows: Let $D$ be
a central division algebra over $k$. Let $D^1$ = the group of reduced
norm 1 elements in $D$. Let $S_0 = \{v \in V\backslash \infty | D
\otimes_k k_v = D_v$ is a division algebra$\}$. Let $D^1 \hookrightarrow
\prod_{v \in S_0} D^*_v$ be the diagonal imbedding. The (locally
compact) topology on $\prod_{v \in S_0} \ D^*_v$ induces a topology on
$D^1$. Then any normal subgroup of $D^1$ is either central (and finite)
or is open in the above topology. (This is the result of the work of
Platonov, Rapinchuk, Margulis, Segev, Seitz and Raghunathan). In the
$k$-rank $\geq 1$ situation the presence of unipotent subgroups holds the
key to the problem. In the case of groups of type $B_n, C_n, D_n$ one
exploits the presence of {\it reflections} in these\break groups.

The cohomology computations were carried out in classical cases by using
the work of Moore. The general situation needs some more refined
machinery -- the Bruhat--Tits theory of buildings associated to groups
over local fields.

If one knows the expectation to hold for an $S$, it will hold for larger
$S$. So the aim would be to handle finite $S$.

The techniques used to handle the case $S=V$ can be used to handle some
kinds of $S$: for example if $K$ is a finite extension and $S = \{v \in
V | K$ does not split completely in $v\}$, then $C(S, G) = 1$ if $G$
is of inner type $A_n$ (and $k$-rank $G= 0$).

When $S$-rank $G = 0$, any $S$-arithmetic group is finite and $C(S, G)$
is trivial.

In the case of $S$-rank $G=1$ some partial results are known. One
expects that $C(S, G)$ in this case is infinite. And this has been shown
to hold in many (classical) cases. One method is to use the following
result: if $C(S, G)$ is finite then $\Gamma^{ab}$ is finite for any
$S$-arithmetic $\Gamma$. One exhibits then $S$-arithmetic $\Gamma$ in
certain $G$ with $\Gamma^{ab}$ infinite thereby showing $C(S, G)$ is not
finite.

Going back to cohomology computations, one has a good understanding of
the group $\widehat{G}_{c,S}(k)$. It is a `restricted direct product' of
$G(k_v), v \not\in S$. Here $G(k_v)$ is the group of the $k_v$-points,
$k_v$ being the (locally compact) completion of $k$ at $v$ and $G(k_v)$
is the group of $k_v$-points of $G$ equipped with its natural locally
compact topology. This reduces the computations to that of $H^2 (G(k_v),
\QQ/\ZZ)$. Moore carried out the computations in the case when $G$ is
split and Deodhar when $G$ is quasi-split. For dealing with the general
case one uses the Bruhat--Tits buildings: These are {\it contractible}
simplicial complexes on which $G(k_v)$ acts. One compares the cohomology
of $G(k_v)$ with that of an imbedded quasi-split subgroup $H(k_v) \subset
G(k_v)$.

Evidently this gets too technical to interest a general audience.

\section*{References}

I give below a fairly comprehensive list of references dealing with the
congruence subgroup problem. In the main body of the paper the precise
references are not given -- only names of some authors are mentioned.
References 59 and 63 below are detailed surveys.

\begin{enumerate}
\renewcommand\labelenumi{[\arabic{enumi}]}
\leftskip .5pc
\item Bak~A, Le probl\'eme des sous-groupes de congruences et le
probl\'eme metaplectique pour les groupes classiques de rang $>1$, {\it
C.R. Acad. Sci. Paris} {\bf 292} (1981) 307--310

\item Bak~A and Rehmann~U, The congruence subgroup and metaplectic
problems for $SL_{n\geq 2}$ of division algebras, {\it J. Algebra} {\bf
78} (1982) 475--547

\item Bass~H, Lazard~M and Serre~J-P, Sous-groupes d'indices finis dans
$SL(n,\ZZ)$, {\it Bull. Am. Math. Soc.} {\bf 70} (1964) 385--392

\item Bass~H, Milnor~J and Serre~J-P, Solution of the congruence
subgroup problem for $SL_{n}(n\geq 3)$ and $Sp_{2n}(n\geq 2)$, {\it Pub.
Math. IHES} {\bf 33} (1967) 59--137; see also {\bf 44} (1974) 241--244

\item Borel~A, Introduction aux groupes arithm\'etiques (Paris: Hermann)
(1969)

\item Borel~A and Harish-Chandra, Arithmetic subgroups of algebraic
groups, {\it Ann. Math.} {\bf 75} (1962) 485--535

\item Borel~A and Tits~J, Groupes R\'eductifs, {\it Publ. Math. IHES}
{\bf 27} (1965) 55--110

\item Borovoi~M~V, Abstract simplicity of some simple anisotropic
algebraic groups over number fields, {\it Sov. Math. Dokl.} {\bf 32}
(1985) 191--193

\item Birtto~J, On defining a subgroup of the special linear group by a
congruence, {\it J. Ind. Math. Soc.} {\bf 40} (1976) 235--243

\item Bruhat~F and Tits~J, Groups r\'eductifs sur un corps local, I:
Donn\'ees radicielles valu\'e\'es, {\it Publ. Math. IHES} {\bf 41}
(1972) 5--251; II: Sch\'emas en groupes. Existence d'une donn\'ee
radicielles, {\it Publ. Math. IHES} {\bf 60} (1984) 5--184

\item \^Cernusov~V, On the projective simplicity of certain groups of
rational points over algebraic number fields, {\it Math. USSR Izv.} {\bf
34} (1990) 409--423

\item Chevalley~C, Sur certains groupes simples, {\it Tohoku J. Math.}
{\bf 7(2)} (1955) 14--62

\item Chevalley~C, Deux th\'eor\`emes d'arithm\'etiques, {\it J. Math.
Soc. Japan} {\bf 3} (1951) 36--44

\item Corlette~K, Archimedean super-rigidity and hyperbolic geoemetry,
{\it Ann. Math.} {\bf 135} (1992) 165--182

\item Delaroche~C and Kirillov~A, Sur les relations entrel'espace dual
d'un groupe et la structure de ses sous-groupes ferm\'es, Expos\'e 343,
S\'eminaire Bourbaki (1967--68)

\item Deligne~P, Extensions centrales non r\'esiduellement finis de
groupes arithm\'etiques, {\it C.R. Acad. Sci. Paris, Ser A-B} {\bf
287(4)} (1978) A203--208

\item Deodhar~V~V, On central extensions of rational points of algebraic
groups, {\it Am. J. Math.} {\bf 100} (1978) 303--386

\item Dickson~L~E, Linear group (Leipzig: Teubner) (1901)

\item Dieudonn\'e~J, La g\'eometrie des groupes classiques (Berlin:
Springer-Verlag) (1955)

\item Fricke~R and Klein~F, Vorlesungen \"uber die Theorie der
automorphen Funktionen, Band I: Die gruppentheoretischen Grundlagen,
Band II: Die funktionentheoretischen Ausf\"{u}hrungen und die Andwendungen
(German) (Stuttgart: B.G. Teubner Verlagsgesellschaft) (1965)

\item Harder~G, Minkowskische Reduktionstheorie \"uber
Funktionenk\"orpern, {\it Invent Math.} {\bf 7} (1969) 33--54

\item Kazhdan~D~A, On the connection between the dual space of the group
with the structure of the closed subgroups, {\it Funct. Anal. Appl.}
{\bf 1} (1967) 71--74 (Russian)

\item Kazhdan~D~A, Some application of the Weil representation, {\it J.
Analyse Mat.} {\bf 32} (1977) 235--248

\item Kazhdan~A~A and Bernstein~I~N, The one-dimensional cohomology of
discrete subgroups (Russian), {\it Funkcional. Anal. i Pril.} {\bf 4}
(1970) 1--5

\item Kneser~M, Orthogonale Gruppen \"uber algebraischen Zahlk\"orpern,
{\it Crelles J.} {\bf 196} (1956) 213--220

\item Kneser~M, Normalteiler ganzzahliger Spingruppen, {\it Crelles J.}
{\bf 311/312} (1979) 191--214

\item Kostant~B, On the existence and irreducibility of certain series
of representation, {\it Bull. Am. Math. Soc.} {\bf 75} (1969) 627--642

\item Li~J-S, Non-vanishing theorems for the cohomology of some
arithmetic quotients, {\it J. Reine Angew. Math.} {\bf 428} (1992) 111--217

\item Lubotzky~A, Group presentation, $p$-adic analytic groups and
lattices in $SL(2, \CC)$, {\it Ann. Math.} {\bf 118} (1983) 115--130

\item Margulis~G~A, Arithmetic properties of discrete groups, {\it
Russian Math. Surveys} {\bf 29} (1974) 107--165

\item Margulis~G~A, Arithmeticity of non-uniform lattices in weakly
non-compact groups (Russian), {\it Funkcional Anal. i Prilozen} {\bf 9}
(1975) 35--44

\item Margulis~G~A, Arithmeticity of the irreducible lattices in the
semisimple groups of rank greater than 1, {\it Inv. Math.} {\bf 76}
(1984) 93--120

\item Margulis~G~A, Finiteness of quotient groups of discrete
subgroups, {\it Funct. Anal. Appl.} {\bf 13} (1979) 178--187

\item Margulis~G~A, On the multiplicative group of a quaternion algebra
over a global field, {\it Soviet Math. Dokl.} {\bf 21} (1980) 780--784

\item Matsumoto~H, Sur les groupes arithm\'etiques des groupes
semisimples d\'eploy\'es, {\it Ann. Sci. Ecole Norm. Sup. (4$^e$ ser.)}
{\bf 2} (1969) 1--62

\item Mennicke~J, Finite factor groups of the unimodular groups, {\it
Ann. Math.} {\bf 81} (1965) 31--37

\item Mennicke~J, Zur Theorie der Siegelschen Modulgruppe, {\it Math.
Annalen.} {\bf 159} (1965) 115--129

\item Millson~J, Real vector bundles with discrete structure group, {\it
Topology} {\bf 18} (1979) 83--89

\item Millson~J, On the first Betti number of constant negatively curved
manifolds, {\it Ann. Math.} {\bf 104} (1976) 235--247

\item Moore~C~C, Extensions and low dimensional cohomology theory of
locally compact groups I and II, {\it Trans. Am. Math. Soc.} {\bf 113}
(1964) 40--86; III, {\it Trans. Am. Math. Soc.} {\bf 221} (1976) 1--38

\item Moore~C~C, Group extensions of $p$-adic and adelic groups, {\it
Publ. Math. IHES} {\bf 35} (1969) 5--70

\item Platonov~V~P, The problem of strong approximation and the
Kneser--Tits conjecture for algebraic groups, {\it Math. USSR Izv.} {\bf
3} (1969) 1139--1147

\item Platonov~V~P, Arithmetic and structure problems in linear
algebraic groups, {\it Proc. ICM Vancouver} {\bf 1} (1974) 471--476

\item Platonov~V~P, The Tannaka--Artin problem, {\it Sov. Math. Dokl.}
{\bf 16} (1975) 468--473

\item Platonov~V~P and Rapinchuk~A~S, Algebraic groups and number theory
(Academic Press) (1991)

\item Platonov~V~P and Rapinchuk~A~S, On the group of rational points of
three dimensional groups, {\it Sov. Math. Dokl.} {\bf 20} (1979)
693--697

\item Platonov~V~P and Rapinchuk~A~S, The multiplicative structure of
division algebras over number fields and the Hasse norm principle, {\it
Sov. Math. Dokl.} {\bf 26} (1982) 388--390

\item Prasad~G, Strong approximation, {\it Ann. Math.} {\bf 105} (1977)
553--572

\item Prasad~G, A variant of a theorem of Calvin Moore, {\it C.R. Acad.
Sci. (Paris) Ser I} {\bf 302} (1982) 405--408

\item Prasad~G and Raghunathan~M~S, Topological central extensions of
semi-simple groups over local fields, {\it Ann. Math.} {\bf 119} (1984)
143--268

\item Prasad~G and Raghunathan~M~S, On the congruence subgroup problem:
Determination of the `metaplectic kernel', {\it Inv. Math.} {\bf 71}
(1983) 21--42

\item Prasad~G and Raghunathan~M~S, On the Kneser--Tits problem, {\it
Comm. Math. Helv.} {\bf 60} (1985) 107--121

\item Prasad~G and Rapinchuk~A~S, Computation of the metaplectic kernel,
{\it Publ. Math. IHES} {\bf 84} (1997) 91--187

\item Raghunathan~M~S, On the congruence subgroup problem I, {\it Publ.
Math. IHES} {\bf 46} (1976) 107--161

\item Raghunathan~M~S, On the congruence subgroup problem II, {\it Inv.
Math.} {\bf 85} (1986) 73--117

\item Raghunathan~M~S, Torsion in co-compact lattices of spin $(2, N)$,
{\it Math. Ann.} {\bf 266} (1984) 403--419

\item Raghunathan~M~S, On the group of norm 1 elements in a division
algebra, {\it Math. Ann.} {\bf 279} (1988) 457--484

\item Raghunathan~M~S, A note on generators for arithmetic subgroups of
algebraic groups, {\it Pac. J. Math.} {\bf 152} (1992) 365--373

\item Raghunathan~M~S, The congruence subgroup problem, Proceedings of
the Hyderabad Conference on Algebraic Groups, Hyderabad, India, Dec.
1989, Manoj Prakashan, India

\item Rapinchuk~A~S, On the congruence subgroup problem for algebraic
groups, {\it Dokl. Akad. Nauk. SSSR} {\bf 306} (1989) 1304--1307

\item Rapinchuk~A~S, Multiplicative arithmetic of division algebras over
number fields and the metaplectic problem, {\it Math. USSR Izv.} {\bf
31} (1988) 349--379

\item Rapinchuk~A~S, Combinatorial theory of arithmetic groups, Preprint
(Acad. of Sciences BSSR) (1990)

\item Rapinchuk~A~S, Congruence subgroup problem for algebraic groups,
old and new, {\it Journe\'es Arithm\'etiques} (Geneva) (1991), {\it
Ast\'erisque} {\bf 209(1)} (1992) 73--84

\item Rapinchuk~A~S, Segev~Y and Seitz~G~M, Finite quotients of the
multiplicative group of a finite dimensional division algebra are
solvable, {\it J. Am. Math. Soc.} {\bf 15(4)} (2002) 929--978
(electronic)

\item Rapinchuk~A~S and Segev~Y, Valuation-like maps and the congruence
subgroup property, {\it Invent. Math.} {\bf 144(3)} (2001) 571--607

\item Segal~D, Congruence topologies in commutative rings, {\it Bull.
London Math. Soc.} {\bf 11} (1979) 186--190

\item Segev~Y, On finite homomorphic images of the multiplicative group
of a division algebra. {\it Ann. Math. (2)} {\bf 149(1)} (1999) 219--251

\item Serre~J-P, Le proble\'me des groupes de congruence pour $SL_2$,
{\it Ann. Math.} {\bf 92} (1970) 489--527

\item Serre~J-P, Sur les groupes de congruence des varie\'te\'s
abelienne, {\it Izv. Akad. Nauk. SSSR} {\bf 28} (1964) 3--18; II, {\it
Akad. Nauk. SSSR} {\bf 35} (1971) 731--735

\item Serre~J-P, Le proble\'me des groupes de congruence pour $SL_2$,
{\it Ann. Math.} {\bf 92} (1970) 489--527

\item Serre~J-P, Sur les groupes de congruence des varie\'te\'s
abelienne, {\it Izv. Akad Nauk. SSSR} {\bf 28} (1964) 3--18; II, {\it
Akad. Nauk. SSSR} {\bf 35} (1971) 731--735

\item Steinberg~R, Ge\'ne\'erateurs, re\'flations et reve\^tements de
groupes alge\'briques, Colloque de Bruxelles, CNRB (1962) 113--127

\item Sury~B, Congruence subgroup problem for anisotropic groups over
semilocal rings, {\it Proc. Indian Acad. Sci. (Math. Sci.)} {\bf 101}
(1991) 87--110

\item Swan~R~G, Generators and relations for certain special linear
groups, {\it Adv. Math.} {\bf 6} (1971) 1--77

\item Tits~J, Algebraic and abstract simple groups, {\it Ann. Math.}
{\bf 80} (1964) 313--329

\item Tits~J, Classification of algebraic semi-simple groups, {\it Proc.
Symp. Pure Math.} {\bf 9} (1966) (AMS)

\item Tits~J, Systems ge\'nerateurs de groupes de congruence, {\it C.R.
Acad. Sci. Paris Series A} {\bf 283} (1976) 693--695

\item Toledo~D, Projective varieties with non-residually finite
fundamental group, {\it Inst. Hautes \'{E}tudes Sci. Publ. Math.} {\bf
77} (1993) 103--119

\item Tomanov~G, On the congruence subgroup problem for some anisotropic
algebraic groups over number fields, {\it Crelles J.} {\bf 402} (1989)
138--152

\item Tomanov~G, Projective simplicity of groups of rational points of
simply connected algebraic groups defined over number fields, Topics in
Algebra Part~2 (Warsaw) (1998) 455--466, (Banach Centre Publications)
26, Part~2, PWN, Warsaw (1990)

\item Tomanov~G, On the group of elements of reduced norm $l$ in a
division algebra over a global field (Russian), {\it Izv. Akad. Nauk.
SSSR Ser. Mat.} {\bf 55} (1991) 917; translation in {\it Math.
USSR-Izv.} {\bf 39} (1992) 895--904

\item Vaserstein~L, The structure of classical arithmetic group of rank
greater than one, {\it Math. USSR Sb.} {\bf 20} (1973) 465--492

\item Venkataramana~T~N, On super-rigidity and arithmeticity of lattices
in semisimple groups, {\it Invent Math.} {\bf 92} (1988) 255--306

\item Venkataramana~T~N, On systems of generators of arithmetic
subgroups of higher rank groups, {\it Pac. J. Math.} {\bf 166} (1994)
193--212

\item Wall~C~T~C, On the commutator subgroups of certain unitary groups,
{\it J. Algebra} {\bf 27} (1973) 306--310

\item Wall~G~E, The structure of a unitary factor group, {\it Publ.
Math. IHES} {\bf 1} (1959)

\item Wang~S~P, The dual space of a semisimple Lie group, {\it Am. J.
Math.} {\bf 91} (1969) 921--937
\end{enumerate}

\end{document}